\documentclass[12pt,bezier]{amsart}
\usepackage{amsfonts,latexsym,amscd}
\newtheorem{theorem}{Theorem}[section]

\newtheorem{proposition}[theorem]{Proposition}

\newtheorem{lemma}[theorem]{Lemma}

\newcommand{\thm}[1]{Theorem \ref{#1}}
\newcommand{\eqn}[1]{Equation \ref{#1}}
\newcommand{\prop}[1]{Proposition \ref{#1}}
\newcommand{\lma}[1]{Lemma \ref{#1}}

\newcommand{\Frac}{\displaystyle\frac}
\newcommand{\rarrow}{\rightarrow}
\newcommand{\fd}{\operatorname{fd}}

\newcommand{\spec}{\operatorname{Spec}}
\newcommand{\conn}{\operatorname{conn}}

\newcommand{\Ker}{\operatorname{Ker}}

\newcommand{\wt}{\operatorname{wt}}
\newcommand{\Hom}{\operatorname{Hom}}
\newcommand{\chr}{\operatorname{char}}
\newcommand{\depth}{\operatorname{depth}}
\newcommand{\tor}{\operatorname{Tor}}
\newcommand{\aqdim}{\operatorname{AQ-dim}}
\newcommand{\la}{\operatorname{\!\langle\!}}
\newcommand{\ra}{\operatorname{\!\rangle}}
\newcommand{\rb}{\operatorname{\!\rbrack}}
\newcommand{\lb}{\operatorname{\!\lbrack\!}}

\newcommand{\F}{{\mathbb F}}

\newcommand{\fm}{\mathfrak{m}}
\newcommand{\fn}{\mathfrak{n}}

\newcommand{\calS}{{\mathcal S}}

\newcommand{\calB}{{\mathcal B}}

\begin{document}

\title[Nilpotency in the Homotopy of Simplicial Algebras]
{Nilpotency in the Homotopy of Simplicial Commutative Algebras}
\author{James M. Turner}
\address{Department of Mathematics\\
Calvin College\\
3201 Burton Street, S.E.\\
Grand Rapids, MI 49546}
\email{jturner@calvin.edu}
\thanks{Partially supported by a grant from the National Science 
Foundation (USA)}
\date{\today}

\begin{abstract}
In this paper, we continue a study of simplicial commutative algebras 
with finite Andr\'e-Quillen homology, that was begun in \cite{Tur}. 
Here we restrict our focus to simplicial algebras having characteristic 2.
Our aim is to find a generalization of the main theorem in \cite{Tur}. 
In particular, we replace the finiteness condition on homotopy with a 
weaker condition expressed in terms of nilpotency for the action of the 
homotopy operations. Coupled with the finiteness assumption on
Andr\'e-Quillen homology, this nilpotency condition provides a way to 
bound the height at which the homology vanishes. As a 
consequence, we establish a special case of an open conjecture of Quillen.
\end{abstract}

\maketitle

\section*{Introduction}

Throughout this paper, unless otherwise stated, all rings and algebras 
are commutative.

Given a simplicial supplemented $\ell$-algebra $A$, with $\ell$ a 
field having non-zero characteristic, it was shown, in \cite{Tur}, 
that if its total Andr\'e-Quillen homology $H^{Q}_{*}(A)$ is 
{\it finite} (as a graded $\ell$-module) then its homotopy $\pi_{*}A$ being 
finite as well implies that $H^{Q}_{*}(A)$ is concentrated in degree 1.
In this paper, we seek to find a generalization of this result by 
weakening the finiteness condition on homotopy. Thus we need to focus 
more on its internal structure. As such we restrict our attention to 
the case where $\chr \ell = 2$ in order to take advantage the rich 
theory available in \cite{Dwy} and \cite{Goe}.

To be more specific, given such a simplicial algebra $A$ having 
characteristic 2, M. Andr\'e \cite{And2} showed that the homotopy groups 
$\pi_{*}A$ have the structure of a divided 
power algebra. Furthermore, W. Dwyer \cite{Dwy} showed that there are 
natural maps
$$
\delta_{i}: \pi_{m}A \to \pi_{m+i}A, \hspace{0.5in} 2\leq i\leq m,
$$
which are homomorphisms for $i<m$ and $\delta_{m} = \gamma_{2}$, the 
divided square. All resulting primary operations can now be described 
in terms of linear combinations of composites of the 
$\delta_{i}$s. Furthermore, there are Adem relations which allow any 
such composite to be described in terms of admissible composites. This 
gives the set $\calB$ of operations a non-commutative ring 
structure. $\pi_{*}A$ then becomes an algebroid over $\calB$. Moreover, 
the module of indecomposables $Q\pi_{*}A$ inherits the structure of 
an unstable $\calB$-module.

A different perspective of Dwyer's operations can be taken in the 
following way. For $0\leq i \leq n-2$ define
$$
\alpha_{i}: \pi_{n}A \to \pi_{2n-i}A
$$
by $\alpha_{i}(x) = \delta_{n-i}(x)$. So, for example, 
$\alpha_{0}(x) = \gamma_{2}(x)$. Written in this way, 
iteration of these reindexed Dwyer operations need not be 
nilpotent. Nevertheless, the main theorem of this paper shows
a connection between the vanishing of Andr\'e-Quillen homology and 
the nilpotency of the Dwyer operations.

\bigskip
\noindent {\bf Theorem A:} {\it \, Let $\ell$ be a field of 
characteristic 2 and let $A$ be a connected simplicial 
supplemented $\ell$-algebra such that the total 
Andr\'e-Quillen homology $H^{Q}_{*}(A)$ is finite. 
Then a nilpotent action of $\alpha_{n-2}$ on $Q\pi_{*}A$ implies that
$H^{Q}_{s}(A) = 0$ for $s\geq n$.
}
\bigskip

As a consequence, the following strengthens the main theorem 
of \cite{Tur} and resolves a conjecture posed in \cite[4.7]{Tur2} at 
the prime 2:

\bigskip 

\noindent {\bf Corollary:} {\it \, Let $A$ be as in Theorem A and
suppose that $\pi_{*}A$ is locally nilpotent as a divided power algebra. 
Then $H^{Q}_{s}(A) = 0$ for $s\geq 2$. 
}
\bigskip 

\noindent {\bf Note:} The restriction to characteristic 2 is due to 
the need for an Adem relations among the homotopy operations 
$\delta_{i}$ which insures that arbitrary composites can be written in 
terms of admissible operations. At the prime 2, this was established 
in the work of Dwyer \cite{Dwy} and Goerss-Lada \cite{GL}. At odd 
primes, Bousfield's work \cite{Bou} gives a preliminary version of 
Adem relations, but a final version still awaits to be produced.

\bigskip 

\subsection*{Connection to conjectures of Quillen}

For a simplicial algebra $A$ over a ring $R$ 
the {\it Andr\'e-Quillen homology}, $D_{*}(A|R;M)$, of $A$ over $R$ with 
coefficients in an $A$-module $M$ was first defined by M. Andr\'e and 
D. Quillen \cite{And, Qui2, Qui3}.  
In particular, for a simplicial supplemented $\ell$-algebra $A$, we 
write
$$
H^{Q}_{*}(A) := D_{*}(A|\ell;\ell).
$$

Next, recall that a homomorphism $\varphi: R \to S$ of Noetherian rings is 
{\it essentially of finite type} if for each $\fn \in \spec S$ there 
is a factorization
\begin{equation}\label{factor}
R \stackrel{\tau}{\rarrow} R[X]_{\mathfrak N} 
\stackrel{\sigma}{\rarrow} S_{\fn}
\end{equation}
where $R[X] = R[X_{1},\ldots,X_{n}]$ is a polynomial ring, $\mathfrak 
N$ is a prime ideal in $R[X]$ lying over $\fn$, the homomorphism 
$\tau: R \to R[X]_{\mathfrak N}$ is the localization map, and the 
homomorphism $\sigma$ is surjective. Furthermore, we call such a 
homomorphism a {\it locally complete intersection} if, for each $\fn 
\in \spec S$, $\Ker (\sigma)$ is generated by a regular sequence.

In \cite[(5.6, 5.7)]{Qui2}, Quillen formulated the following two 
conjectures on the vanishing of Andr\'e-Quillen homology:

\bigskip
\noindent {\bf Conjecture:} {\it \, Let $\varphi: R \to S$ be a
homomorphism essentially of finite type between Noetherian rings and 
assume further that $D_{s}(S|R;-)  = 0$ for $s\gg 0$. Then:
\begin{enumerate}
	
	\item[I.] $D_{s}(S|R;-) = 0$ for $s\geq 3$;
	
    \item[II.] If, additionally, the flat dimension $\fd_{R} S$ is finite, 
	then $\varphi$ is a locally complete intersection homomorphism.
       	  
\end{enumerate}	  
}
\bigskip

In \cite{Avr}, L. Avramov generalized the notion of local 
complete intersections to arbitrary homomorphisms of Noetherian rings. 
He further proved a generalization of Conjecture II. to such 
homomorphisms. See P. Roberts review \cite{Rob} of this 
paper for an excellent summary of these results and the history behind 
them. A proof of Conjecture II. was also given in \cite{Tur} for 
homomorphisms with target $S$ having non-zero characteristic.

\bigskip 

We now indicate how Theorem A bears on providing a resolution to the 
above Conjecture. To formulate this, let $R \to (S,\ell)$ be a 
homomorphism of local rings with $\chr \ell = 2$. Let 
$$
\vartheta : (Q\tor^{R}_{*}(S,\ell))_{m} \to (Q\tor^{R}_{*}(S,\ell))_{2m-1}
$$
be the operation induced by $\alpha_{1}$. We call this 
operation the {\it Andr\'e operation} since it generalizes the 
operation studied by M. Andr\'e \cite{And3} when $S = \ell$ and $m = 
3$. 

\bigskip


\noindent {\bf Theorem B:} {\it \, Let $\varphi: R\to (S,\ell)$ be a 
surjective homomorphism of local rings with $\chr \ell = 2$ and assume further 
that $D_{s}(S|R;\ell) = 0$ for $s\gg 0$. Then 
\begin{enumerate}
	
	\item $D_{s}(S|R;\ell) = 0$ for $s\geq 3$ if and only if the Andr\'e 
	operation $\vartheta$ acts nilpotently on $Q\tor^{R}_{*}(S,\ell)$;
	
	\item $\varphi$ is a complete intersection if and only if the divided
	square $\gamma_{2}$ acts nilpotently on $Q\tor^{R}_{*}(S,\ell)$.
	
\end{enumerate}	
}
\bigskip 

As an application of Theorem B, the following proves the vanishing 
portions of the conjecture for certain homomorphisms in 
characteristic 2:

\bigskip
\noindent {\bf Theorem C:} {\it \, Let $R\to S$ be a 
homomorphism essentially of finite type between Noetherian rings of
characteristic 2 and assume further that $D_{s}(S|R;-) = 0$ for $s\gg 0$. 
Then:
\begin{enumerate} 
	
	\item $D_{s}(S|R;-) = 0$ for $s\geq 3$ provided $R\to S$  
	is a homomorphism of Cohen-Macaulay rings;
	
	\item If the flat dimension $\fd_{R} S$ is finite, then $\varphi$ is a 
	locally complete intersection.
	
\end{enumerate}	
	
}
\bigskip 

\noindent {\it Note:} L. Avramov showed in \cite{Avr} that Conjecture I. holds 
when either $R$ or $S$ is a locally complete intersection.  
More recently,  L. Avramov and S. Iyengar 
\cite{AI1} have strengthened this special case of Conjecture I.  
by showing that it holds for homomorphisms $R \to S$ for which 
there exists a composite $Q \to R \to S$ which is a local complete 
intersection homomorphism. See \cite{AI2} for a more leisurely 
discussion of their results.

\subsection*{Organization of this paper}
The first section reviews the properties of the homotopy and 
Andr\'e-Quillen homology for simplicial commutative algebras and the 
methods for computing them, particularly in characteristic 2. The 
next section then focuses on a device called the character map 
associated to simplicial algebras having finite Andr\'e-Quillen 
homology. After showing that Dwyer's operations possess certain 
annihilation properties, we show that the character map can be highly 
non-trivial. From this Theorem A easily follows. This enables us, in 
the third section, to prove Theorem B. Finally, the last section 
begins with establishing a chain level criterion for the nilpotency 
of Andr\'e's operation. After a brief excursion into commutative 
algebra, we prove a special case of Theorem C and then show how the
general case follows.

\subsection*{Acknowledgements}
The author would like to thank Lucho Avramov for sharing his
expertise and insights and for a thorough and critical reading of an earlier 
draft of this paper.
He would also like to thank Paul Goerss and Jean Lannes for helpful advice 
and discussions during work on this project. 

\section{Homotopy and homology of simplicial commutative algebras in 
characteristic 2}

\subsection{Operations on chains and homotopy}

Again let $A$ be a simplicial algebra of characteristic 2. We 
describe the algebra structure associated to $A$ at two levels: on the 
associated chain complex and on the associated homotopy groups. While only 
the latter will be needed in the process of proving Theorem A, the 
former description will be needed in the subsequent applications. 

Let $V$ be a simplicial vector space, over the field $\F_{2}$,
and let $C(V)$ denote its associated chain complex. 
In \cite{Dwy}, Dwyer constructs natural chain maps
\begin{equation}\label{d1}
\Delta^{k}: (C(V)\otimes C(V))_{n+k}\to C(V\otimes V)_{n}
\end{equation}
for all $0\leq k\leq n$. They satisfy the relations
\begin{equation}\label{rel1}
\Delta^{0}+T\Delta^{0}T+\phi_{0} = \Delta
\end{equation}
and
\begin{equation}\label{rel2}
\Delta^{k}+T\Delta^{k}T+\phi_{k} = \partial \Delta^{k-1}+\Delta^{k-1}\partial. 
\end{equation}
Here $T:C(V)\otimes C(V) \to C(V)\otimes C(V)$ and $T : C(V\otimes 
V)\to C(V\otimes V)$ are the twist maps. Also, $$\phi_{k}: C(V)\otimes 
C(V) \to C(V\otimes V), \quad k\geq 0$$ is the degree $-k$ map that 
is zero on $[C(V)\otimes C(V)]_{m}$ for $m\neq 2k,$ and, in degree 
$2k$, is the projection on one factor:
$$
[C(V)\otimes C(V)]_{2k} = \oplus_{p+q=2k}V_{p}\otimes V_{q} \to 
V_{k}\otimes V_{k}.
$$
Finally $\Delta: C(V)\otimes C(V) \to C(V\otimes 
V)$ is the Eilenberg-Zilber map.

Now given a simplicial $\F_{2}$-algebra $(A,\mu)$, define the maps
\begin{equation}\label{alf}
\alpha_{i} : C(A)_{n} \to C(A)_{2n-i} {\it \quad for \quad} 0\leq i \leq n-1
\end{equation}
by
\begin{equation}\label{alfdef}
x \to \mu \Delta^{i}(x\otimes x)+\mu \Delta^{i-1}(x\otimes \partial x).
\end{equation}

To describe the algebra structure on the associated chains of $A$, 
recall that a {\it dg $\Gamma$-algebra} is a non-negatively 
differentially graded $\F_{2}$-algebra $(\Lambda,\partial)$ together with maps
$$
\gamma_{k} : \Lambda_{n} \to \Lambda_{kn} \quad {\mbox for} \quad 
k\geq 0 \quad {\mbox and} \quad n\geq 2,
$$
satisfying the following relations

\begin{enumerate}
	
	\item $\gamma_{0}(x) = 1$ and $\gamma_{1}(x) = x$
	
	\item $\gamma_{h}(x)\gamma_{k}(x) = \binom{h+k}{h}\gamma_{h+k}(x)$
	
	\item $\gamma_{k}(x+y) = \sum_{r+s = k}\gamma_{r}(x)\gamma_{s}(x)$
	
	\item $\gamma_{k}(xy) = 0$ for $k\geq 2$ and $x,y \in \Lambda_{\geq 
	1}$
	
	\item $\gamma_{k}(xy) = x^{k}\gamma_{k}(y)$ for $x\in \Lambda_{0}$ 
	and $y\in \Lambda_{\geq 2}$
	
	\item $\gamma_{k}(\gamma_{2}(x)) = \gamma_{2k}(x)$
	
	\item $\partial\gamma_{k}(x) = (\partial x)\gamma_{k-1}(x)$.
	
\end{enumerate}	

\begin{proposition}\label{dwylmma}
	Let $A$ be a simplicial $\F_{2}$-algebra. 
	\begin{enumerate}
		\item The chain complex $C(A)$ possesses a dg $\Gamma$-algebra 
		with $\alpha_{0} = \gamma_{2}$;
		
		 \item For $x\in C(A)_{n}$ and $0<i<n-1$, $$\partial(\alpha_{i}(x)) = 
         \alpha_{i-1}(\partial (x)), {\it \quad for \quad} 0<i<n-1,$$ 
         $$\partial(\alpha_{0}(x)) = x\cdot \partial (x), {\it \quad and \quad} 
		 \partial(\alpha_{n-1}(x)) = \alpha_{n-2}(\partial x) + x^{2}.$$
	\end{enumerate}	 
\end{proposition}

\noindent {\it Proof:} 1. See \cite{And2} and \cite[\S 2 and 3]{Goe}.

2. By \ref{rel2}, \ref{alfdef}, and the Leibniz 
rule for $\partial$, 
\begin{eqnarray*}
\partial\alpha_{i}(x) & = & \mu \partial \Delta^{i}(x\otimes x)+\mu \partial 
                            \Delta^{i-1}(x\otimes \partial x) \\
					  &	= & \mu \Delta^{i}\partial (x\otimes x)+\mu 
					        \Delta^{i-1}\partial (x\otimes \partial x) +
							\mu \Delta^{i}(\partial x\otimes x)+\mu 
							\Delta^{i}(x\otimes \partial x) \\
					  & = & \mu \Delta^{i-1}(\partial x\otimes \partial x) \\
					  & = & \alpha_{i-1}(\partial x).					  
\end{eqnarray*}
The calculation of $\partial \alpha_{0}(x)$ and $\partial \alpha_{n-1}(x)$ 
from (\ref{rel1}) and (\ref{rel2}) is similar. \hfill $\Box$ 

\bigskip 

\noindent {\bf Note:} For a divided power algebra $\Lambda$ of 
characteristic 2, it is enough to specify the action of divided 
square $\gamma_{2}$ to determine all divided powers. Specifically, for 
$x\in \Lambda$
$$
\gamma_{k} = \gamma_{2}^{s_{1}}(x)\cdot \gamma_{2}^{s_{2}}(x)\ldots 
\gamma_{2}^{s_{r}}(x), 
$$
where $k = 2^{s_{1}}+\ldots +2^{s_{r}}$. Cf. \cite{And2} and \cite[\S 2]{Goe}.

\bigskip 

Now, the $\alpha_{i}$ induce the homotopy operations
$$
\alpha_{i}: \pi_{n}A \to \pi_{2n-i}A, \quad 0\leq i \leq n-2.
$$
Furthermore, Dwyer's {\it higher divided squares} 
$$
\delta_{i}:\pi_{n}A \to \pi_{n+i}A, \quad 2\leq i\leq n
$$
are now defined as
$$
\delta_{i}[x] = [\alpha_{n-i}(x)].
$$
In particular,
$$
\delta_{n}[x] = [\alpha_{0}(x)] = \gamma_{2}[x].
$$

We now summarize the properties of the higher divided squares, as 
established in \cite{Dwy,GL} (see also \cite[\S 2 and 3]{Goe}).

\begin{proposition}\label{dwyprop}
The higher divided squares possess the following properties:
\begin{enumerate}
	
	\item Adem relations:
	
	\begin{enumerate}
		
		\item For $i<2j$, 
           $$\delta_{i}\delta_{j} = \sum_{\frac{i+1}{2}\leq s\leq 
           \frac{i+j}{3}}\binom{j-i+s-1}{j-s}\delta_{i+j-s}\delta_{s}$$
		   
        \item For $j<i$,
		   $$\alpha_{i}\alpha_{j} = \sum_{\frac{i+2j}{3}\leq s\leq 
             \frac{i+j-1}{2}}\binom{i-s-1}{s-j}\alpha_{i+2j-2s}\alpha_{s}$$
			 
	\end{enumerate}		 

    \item Cartan formula: for $x,y \in \pi_{*}A$, 
      $$\delta_{i}(xy) =
             \begin{cases}
              x^{2}\delta_{i}(y) & \quad |x| = 0; \\
              y^{2}\delta_{i}(x) & \quad |y| = 0; \\
              0 & \quad |x|>0,|y|>0.
             \end{cases}$$
			 
\end{enumerate}			 
\end{proposition}

Let $I = (i_{1},\ldots,i_{s})$ be a sequence of positive integers. 
Then call $I$ {\it admissible} provided $i_{t}\geq 2i_{t+1}$ for all 
$1\leq t < s$. Furthermore, define for $I$ its {\it excess} to be 
the integer
$$
e(I) = (i_{1}-2i_{2})+(i_{2}-2i_{3})+\ldots+(i_{s-1}-2i_{s})+i_{s} = 
i_{1}-i_{2}-\ldots-i_{s},
$$
its {\it length} to be the integer $\lambda(I) = s$, and its {\it 
degree} to be the integer
$$
d(I) = i_{1}+\ldots+i_{s}.
$$
Also write
$$
\delta_{I} = \delta_{i_{1}}\ldots\delta_{i_{s}} \qquad \alpha^{I} = 
\alpha_{1}^{i_{1}}\ldots\alpha_{s}^{i_{s}}.
$$

As an application of the Adem relations, given any sequence $I$, $\delta_{I}$ 
can be written as a sum of $\delta_{J}$'s, with each $J$ an admissible 
sequence. Similarly, by another application of the Adem relations, 
any $\alpha_{I}$ can be written as a sum of $\alpha^{J}$'s, with $J$ 
not necessarily admissible.

Finally, denote by $\calB$ the algebra spanned by $\{\delta_{I}|I \; 
\mbox{admissible}\}$. A $\calB$-module $M$ is then called {\it unstable} 
provided $\delta_{I}x = 0$ for any $x\in M_{n}$, whenever $I$ is 
admissible with $e(I)>n$. For example, given a simplicial 
supplemented $\ell$-algebra $A$ then $Q\pi_{*}A$, the module of 
indecomposables, is an unstable $\calB$-module.

\subsection{The homotopy of symmetric algebras}

Let $\ell$ be a field having characteristic 2. 
We provide description of the homotopy groups of a very 
important type of simplicial $\ell$-algebra, namely, 
the homotopy of $S_{\ell}(V)$ - the symmetric algebra, or 
free commutative algebra, generated by a simplicial vector space 
$V$ over $\ell$. When $\ell = \F_{2}$, we simply write $S$ for 
$S_{\ell}$.

If $W$ is a vector space, let $K(W,n)$ denote the simplicial vector 
space with homotopy $\pi_{*}K(W,n) \cong W$, concentrated in degree $n$. 
Then write
\begin{eqnarray*}
	S_{\ell}(W,n) & = & S_{\ell}(K(W,n)) \\ 
	S_{\ell}(n)   & = & S_{\ell}(\ell,n)
\end{eqnarray*}	

By a theorem of Dold \cite{Dold}, there is a functor of graded vector 
spaces $\calS_{\ell}$ and a natural isomorphism 
\begin{equation}\label{dold1}
\pi_{*}S_{\ell}(V) \cong \calS_{\ell}(\pi_{*}V).
\end{equation}
Again, we simply write $\calS$ for $\calS_{\ell}$ when $\ell = \F_{2}$.

To describe the functor $\calS_{\ell}$ on graded vector spaces, we first note 
that it commutes with colimits. Thus we need only describe 
$\calS_{\ell}(F_{\ell}(n))$ where $F_{\ell}(n) \cong \ell \la x_{n} \ra$. 
We now recall the description of this functor when $\ell = \F_{2}$.
A proof of the following can be found, for example, in \cite{Dwy}:

\begin{proposition}\label{dold2}
\begin{eqnarray*}
\calS(F(n)) & \cong & \Gamma[\alpha^{I}(x_{n}): \sigma(I)<n-1]\\
& \cong & \Gamma[\delta_{I}(x_{n}): I \; \mbox{admissible}, 
\; e(I) < n].
\end{eqnarray*}
\end{proposition}

Here $\Gamma[-]$ denotes the free divided power algebra functor. Note 
that $Q\calS(W)$ is a free unstable $\calB$-module, for any positively 
graded vector space $W$.

The functor $\calS$ can be further decomposed as $$\calS(-) = 
\bigoplus_{m\geq 0}\calS_{m}(-).$$ We review its description because 
of its importance below. 

First, define the {\it weight} of an element of $\calS(W)$ as 
follows:
$$
\wt (u) = 1, \quad \wt (uv) = \wt (u) + \wt (v), 
$$
$$
\wt (\delta_{i}(u)) = 2\wt (u) \text{\quad for \quad} u,v \in W.
$$

\begin{proposition}\label{dold3}
Let $W$ be a graded vector space. Then $\calS_{m}(W)$ is the 
subspace $\calS(W)$ spanned by elements of weight $m$.
\end{proposition}

\bigskip 

Now, to describe $\calS_{\ell}$, note first that, the uniqueness of 
adjoint funtors, there is a natural isomorphism of $\ell$-algebras
$$
\eta : S_{\ell}(V\otimes_{\F_2}\ell) \to S(V)\otimes_{\F_2} \ell
$$
where $V$ is any $\F_{2}$-vector space. 

\begin{proposition}\label{dold4}
	The natural isomorphism $\eta$ in turn induces a natural isomorphism 
	$$\eta_{*} : \calS_{\ell}((-)\otimes_{\F_2}\ell) 
	\stackrel{\cong}{\rarrow} \calS(-)\otimes_{\F_2}\ell$$
	of functors to the category of $\Gamma$-algebras. Moreover, for each $i,m\geq 
	2$ and graded $\F_{2}$-vector space $W$, there is a commutative diagram
	$$
    \begin{array}{ccc}
    \calS_{\ell}(W\otimes_{\F_2}\ell)_{n} & 
    \stackrel{\delta_{i}}{\rightarrow\!\!\!\!\rightarrow} & 
	\calS_{\ell}(W\otimes_{\F_2}\ell)_{n+i} \\[1mm]
    \downarrow \hspace*{1pt}
    &&
    \hspace*{1pt} \downarrow \\[1mm]
    \calS(W)_{n}\otimes_{\F_2}\ell & \stackrel{\delta_{i}\otimes F}{\rightarrow\!\!\!\!\rightarrow} & 
    \calS(W)_{n+i}\otimes_{\F_2}\ell \\[1mm]
    \end{array}
    $$
	where $F$ denotes the Frobenius map on $\ell$.
\end{proposition}	

\bigskip

\noindent {\it Proof:} \; Since $\pi_{*}S_{\ell}(-) \cong 
\calS_{\ell}(-)$, by Dold's theorem, then 
the first point regarding $\eta_{*}$ follows from a Kunneth theorem 
argument. To prove the second part, it suffices to prove it for the 
graded vector space $F_{\ell}(n) = F(n)\otimes_{\F_2}\ell$, by a standard 
argument utilizing universal examples. 
In this case, the map $\eta$ extends the map $\ell \la 
x_{n} \ra \to \F_{2} \la x_{n} \ra \otimes_{\F_2} \ell$ sending $ax_{n}$ to 
$x_{n}\otimes a$. By naturality of $\eta_{*}$ and the properties of 
Dwyer's operations, we have

\begin{eqnarray*}
\delta_{i}(\eta_{*}(ax_{n})) & = & \delta_{i}(x_{n}\otimes a) =  
\delta_{i}((x_{n}\otimes 1)(1\otimes a)) \\ 
& = & \delta_{i}(x_{n}\otimes 1)(1\otimes a)^{2} 
= (\delta_{i}(x_{n})\otimes 1)(1\otimes a^{2}) \\
& = & \delta_{i}(x_{n})\otimes a^{2}
\end{eqnarray*}

which is the desired result. \hfill $\Box$

\bigskip

It now follows that, for a simplicial $\ell$-vector space $V$, 
the generators and relations for
$\pi_{*}S_{\ell}(V)$ are completely determined by Dwyer's result 
\prop{dwyprop} and Dold's theorem.

\bigskip

\subsection{Andr\'e-Quillen homology and the fundamental spectral sequence}

We now provide a brief review of Andr\'e-Quillen homology, for 
simplicial supplemented $\ell$-algebras, and the main computational 
device for relating homotopy and homology - Quillen's fundamental 
spectral sequence. Our primary source for this material is \cite{Goe}. 
Cf. also \cite{Qui2, Qui3, Mil}.

Let $A$ be a simplicial supplemented $\ell$-algebra. Then the {\it 
Andr\'e-Quillen homology} of $A$ is defined as the graded vector space 
$$H^{Q}_{*}(A) = \pi_{*}QX,$$ 
where $X$ is a cofibrant replacement of $A$, in the 
closed simplicial model structure on simplicial supplemented 
$\ell$-algebras \cite{Mil, Goe}.

Some standard properties of Andr\'e-Quillen homology are summarized 
in the following:

\begin{proposition}\label{aqprop}\cite[\S 4]{Goe}
	Let $A$ be a simplicial supplemented $\ell$-algebra.
	\begin{enumerate}
		\item If $A = S_{\ell}(V)$, for some simplicial vector space $V$, then 
		$H^{Q}_{*}(A) \cong \pi_{*}V$.
		
		\item Let $A \stackrel{f}{\rarrow} B \stackrel{g}{\rarrow} C$ be a 
         cofibration sequence in the homotopy category of simplicial 
         supplemented algebras.  Then
         there is a long exact sequence
         $$
         \begin{array}{l}
         \cdots \rarrow H^Q_{s+1} (C) \stackrel{\partial}{\rarrow} H^Q_s(A)
         \stackrel{H^Q_*(f)}{\rarrow} H^Q_s(B) \\[3mm]
         \hspace*{20pt} \stackrel{H^Q_*(g)}{\rarrow} H^Q_s(C)
         \stackrel{\partial}{\rarrow} H^Q_{s-1}(C) \rarrow \cdots
         \end{array}
         $$

		 \item If $V$ is a vector space and [\quad,\quad] denotes morphisms in 
		 the homotopy category of simplicial supplemented algebras, then the map 
         $$
         [S_{\ell}(V,n),A] \rarrow \Hom(V, I\pi_n A)
         $$
         is an isomorphism. In particular, $\pi_{n}A = [S_{\ell}(n),A]$.
	\end{enumerate}
\end{proposition}

Another important tool we will need is the notion of {\it connected 
envelopes} for a simplicial supplemented algebra $A$. Cf. \cite[\S 2]{Tur}. 
These are defined as a sequence of cofibrations
$$
A = A(0) \stackrel{j_{1}}{\rightarrow} A(1) \stackrel{j_{2}}{\rightarrow} \cdots
\stackrel{j_{n}}{\rightarrow} A(n) \stackrel{j_{n+1}}{\rightarrow} \cdots
$$
with the following properties: 
\begin{itemize}
\item[(1)] For each $n\geq 1$, $A(n)$ is a $n$-connected. 
\item[(2)] For $s > n$,
$$
H^Q_s A(n) \cong H^Q_sA.
$$
\item[(3)] There is a cofibration sequence
$$
S_{\ell}(H^Q_{n} A, n) \stackrel{f_{n}}{\rarrow} A(n-1) \stackrel{j_{n}}{\rarrow} A(n).
$$
\end{itemize}

The following is proved in \cite{Tur}:

\begin{lemma}\label{ceprop}
	If $H^{Q}_{s}(A) = 0$ for $s>n$ then $A(n-1)\cong S_{\ell}(H^{Q}_{n}(A), 
	n)$ in the homotopy category of simplicial supplemented algebras.
\end{lemma}	

Finally, a very important tool for bridging the Andr\'e-Quillen 
homology to the homotopy of a simplicial algebra is provided by the 
{\it fundamental spectral sequence} of Quillen \cite{Qui2,Qui3}. For 
simplicial supplemented $\ell$-algebras in characteristic 2, we will 
need certain properties of this spectral sequence which can be described by combining 
the results of \cite[\S 6]{Goe} with \prop{dold4}. The following summarizes 
those features that we will need.

\begin{proposition} \label{qsseq}
Let $A$ be a simplicial supplemented $\ell$-algebra. Then there is a 
spectral sequence of algebras
$$
E^{1}_{s,t}A = \calS_{s}(H^{Q}_{*}(A))_{t}\otimes_{\F_2} \ell \Longrightarrow \pi_{t}A
$$
with the following properties:
\begin{enumerate}
	
	\item For $\pi_{0}A \cong \ell$, $E^{1}_{s,t}A = 0$ for $s>t$ and, hence, 
	the spectral sequence converges;
	
	\item The differentials act as $d_{r}: E^{r}_{s,t} \to E^{r}_{s+r,t-1}$;
	
	\item The Dwyer operations $$\delta_{i}: E^{r}_{s,t}\to 
	E^{r}_{2s,t+i}, \quad 2\leq i \leq t$$ have indeterminacy 2r-1 and 
	satisfy the following properties: 
	
	   \begin{enumerate}
		   
		   \item Up to determinacy, the Adem relations and Cartan formula holds;
		   
		   \item If $x\in E^{r} A$ and $2\leq i < t$, then $\delta_{i}(x)$ 
		   survives to $E^{2r}A$ and
		   \begin{eqnarray*}
			   d_{2r}\delta_{i}(x) & = & \delta_{i}(d_{r}x)\\
			   d_{r}\delta_{t}(x) & = & xd_{r}x
		   \end{eqnarray*}
		   modulo indeterminacy;
		   
		   \item The operations on $E^{r}A$, $r\geq 2$, are induced by the 
		   operations on $E^{r-1}A$. The operations on $E^{\infty}A$ are 
		   induced by the operations on $E^{r}A$ with $r<\infty$; and
		   
		   \item The operations on $E^{\infty}A$ are also induced by the 
		   operations on $\pi_{*}A$.
		   
	    \end{enumerate}

\end{enumerate}
\end{proposition}

Recall that, if 
$$
B^{q}_{s,t}\subseteq E^{r}_{s,t}, \quad q\geq r
$$
is the vector space of elements that survive to $E^{q}_{s,t}$ but have 
zero residue class in $E^{q}_{s,t}$, then $y\in E^{r}_{s,t}$ is 
defined {\it up to indeterminacy} $q$ if $y$ is a coset 
representitive for a particular element in $E^{r}_{s,t}/B^{q}_{s,t}$.

\bigskip

\section{Proof of Theorem A}

\subsection{The character map for simplicial algebras with finite 
homology}

Fix a connected simplicial supplemented $\ell$-algebra $A$ with 
$H^{Q}_{*}(A)$ finite as a graded $\ell$-module. Define the {\it 
Andr\'e-Quillen dimension} of $A$ \cite{AI2} to be $$\aqdim A = 
\mbox{max}\{s:H^{Q}_{s}(A)\neq 0\}$$ and define the {\it 
connectivity} of $A$ to be $$\conn A = 
\mbox{min}\{s:H^{Q}_{s}(A)\neq 0\}-1.$$ We assume that $\aqdim A\geq 2$.

Let $n = \aqdim A$. Then the (n-1)-connected envelope $A(n-1)$ has the property
that $$A(n-1) \cong S_{\ell}(H^{Q}_{n}(A),n)$$ in the homotopy category. Cf. 
\lma{ceprop}. 
Thus we have a map $A\to S_{\ell}(H^{Q}_{n}(A),n)$ in the homotopy category 
with the property that it is an $H^{Q}_{n}$-isomorphism. 

We now define the {\it character map} of $A$ to be the resulting 
induced map of unstable $\calB$-modules
\begin{equation}\label{char}
	\Phi_{A} : Q\pi_{*}A \to Q\pi_{*}S_{\ell}(H^{Q}_{n}(A),n)
\end{equation}	

The importance of the character map is established by the following:

\begin{theorem}\label{charprop}
Let $A$ be a connected simplicial 
supplemented $\ell$-algebra having finite Andr\'e-Quillen dimension $n$. 
Then the character map $\Phi_{A}$ is non-trivial. Furthermore, 
$y\in Q\pi_{*}A$ can be chosen so that $$\Phi_{A}(y) = 
\alpha_{n-2}^{s}(x),$$ for some non-trivial $x\in H^{Q}_{n}(A)$ 
and some $s>0$. Thus
$\alpha_{t-2}$ acts non-nilpotently on $y$, for all $2\leq t\leq n$.
\end{theorem}

As an immediate consequence, we are now in a position to supply the 
following:

\bigskip

\noindent {\it Proof of Theorem A:} This follows immediately from \thm{charprop} 
since if $\alpha_{n-2}$ acts nilpotently on any $x\in Q\pi_{*}A$, the same 
must also hold for $\Phi(x)$. Hence it follows that $n>\aqdim A$. \hfill $\Box$

\subsection{Annihilation properties among some homotopy operations}

Before we prove \thm{charprop}, we need to pin down specific 
annihilators of elements of $\calB$. To this end,
define, for $s,t\geq 0$, the operation $$\theta(s,t) = 
\delta_{2^{s+t}}\delta_{2^{s+t-1}}\ldots\delta_{2^{t+1}}.$$ 

\begin{proposition}\label{annprop}
Let $J$ be a finite subset of $\{j|j>2^{t}\}$ and 
let $$\xi = \sum_{j\in J}a_{j}\delta_{j}w_{j}$$ be a linear 
combination of admissible operations, with each $a_{j}\in \ell$. 
Then $\theta(s,t)\xi = 0$ for $s\gg 0$.
\end{proposition} 

\noindent {\it Proof:} It is sufficient to prove the result for 
$\xi = \delta_{j}$ with $j>2^{t}$. Write $j = 2^{t}+n$ with $n\geq 1$. Note first that an application of 
the Adem relations shows that, for any $t$, 
$$\delta_{2^{t+1}}\delta_{2^{t}+1}=\delta_{2^{t+1}}\delta_{2^{t}+2}=0.$$ 
We thus assume, by induction, that, for each $t$ and $0<i<n$, there exists 
$s\gg 0$ such that $$\theta(s,t)\delta_{2^{t}+i}= 0.$$ By another 
application of the Adem relations, we have
$$
\delta_{2^{t+1}}\delta_{2^{t}+n} = \sum_{1\leq r\leq 
\frac{n}{3}}\binom{n+r-1}{n-r}\delta_{2^{t+1}+n-r}\delta_{2^{t}+r}.
$$
Notice that, for each such $r$, $2^{t+1}< 2^{t+1}+n-r < 2^{t+1}+n$. Thus, by induction, we 
can find $s\gg 0$ so that
$$
\theta(s,t+1)\big{(}\sum_{1\leq r\leq 
\frac{n}{3}}\binom{n+r-1}{n-r}\delta_{2^{t+1}+n-r}\delta_{2^{t}+r}\big{)} = 0.
$$
We conclude that $$\theta(s,t)\delta_{2^{t}+n}=\theta(s,t+1)\delta_{2^{t+1}}\delta_{2^{t}+n} = 0.$$ 
\hfill $\Box$

\bigskip

\subsection{Proof of \thm{charprop}}

As we noted previously, it is sufficient to show that, for some $t>0$, 
$\alpha_{n-2}^{t}(x) = 
\delta_{2^{t}}\delta_{2^{t-1}}\ldots\delta_{2}(x) \in 
E^{1}_{2^{t}}$ survives non-trivially to $E^{\infty}$ for some $x\in 
H^{Q}_{n}(A)$, where $n = \aqdim A$. 
Such an element will map non-trivially under $\Phi$ with the desired 
properties, by \prop{dold2}.

The strategy is to examine the induced map of spectral sequences 
$$\{E^{r}A\}\to \{E^{r}S(H^{Q}_{n}(A),n)\}$$ which is split surjective 
at $E^{1}$. The goal is to show 
that the image of the splitting on the indecomposables 
contains a non-trivial infinite cycle with the requisite specifications.

Now, assume $n\geq 2$. Then the result holds when 
$n = m$ and $n = m+1$, where $m = \conn A+1$, because in these cases Quillen's 
spectral sequence collapses \cite{Tur} and, hence, $\Phi_{A}$ is a split
surjection. Thus we can now induct on $n-m$. This further reduces 
to an induction on $\dim_{\ell}H^{Q}_{m}(A)$, which is finite by assumption. 

By the Hurewicz theorem \cite[(8.3)]{Goe}, $\pi_{m}A \cong H^{Q}_{m}(A)$. 
By \prop{aqprop}.3, a choice of a 
basis element $y\in H^{Q}_{m}(A)$ is represented by a map $\sigma : S(m) 
\to A$ of 
simplicial supplemented $\ell$-algebras, by \prop{aqprop}.3. Let 
$B$ be the homotopy cofibre (aka mapping cone \cite[(4.5)]{Goe}) 
of $\sigma$ and let $f:A\to B$ be the induced map. Note that 
there is an identity:
\begin{equation}\label{charid}
	\Phi_{A}=\Phi_{B}(Qf_{*})
\end{equation}	
Then 
$\dim_{\ell}H^{Q}_{m}(B) = 
\dim_{\ell}H^{Q}_{m}(A)-1$. By induction, we assume that, for $x\in 
H^{Q}_{n}(A)$, $\theta(b,0) 
x\in E^{1}_{2^{b},t}B$ survives non-trivially to $E^{\infty}B$
and determines an element $y^{\prime}\in \pi_{t}B$ such that
\begin{equation}\label{indstep}
\Phi_{B}(y^{\prime}) = \alpha_{n-2}^{b}(x). 
\end{equation}

Now, in the spectral sequence for $A$, \prop{qsseq}.3 (b) tells us, by 
induction, that $\theta(b,0) x$ survives to some 
$E^{r}_{2^{b},t}A$ with $r\geq 2^{b}$. Furthermore, by \prop{dold3},
the differential satisfies:
$$
d_{r}\left [\theta(b,0) x \right ] =
\begin{cases}
  \left [uy+\xi y \right ] & \quad r = 2^{a}-2^{b}, \; a>b; \\
  \left [ uy\right ] & \quad \text{otherwise},
\end{cases}
$$
with $u \in E^{1}_{2^{b}+r-1}$ and $\xi$ a linear combination of admissible Dwyer 
operations. Our goal is to show that there exists 
$s\gg b$ so that $\theta(s,0) x$ is an infinite 
cycle.  We examine the above cases on $r$ in reverse 
order. 
\bigskip

\noindent ${\bf r\neq 2^{a}-2^{b}:}$ We use induction on $t-2^{b}-r$, by 
first noting that for $t-2^{b}-r=0$, $d_{r}\lb\theta(b,0)x\rb \in 
E^{r}_{2^{b}+r,t-1}=0$. 
Since $\theta(b+1,0)=\delta_{2^{b+1}}\theta(b,0)\neq 0$ then 
$\delta_{2^{b+1}}\theta(b,0) x$ survives to $E^{2r}$ 
by \prop{qsseq}.3 (b). By \prop{qsseq}.1 and \prop{qsseq}.3 (a) and (d),  
$$d_{2r}\lb\delta_{2^{b+1}}\theta(b,0) 
x\rb = \delta_{2^{b+1}}\lb d_{r}\theta(b,0) x\rb = 0.$$ Hence 
$\delta_{2^{b+1}}\theta(b,0) x$ survives
to $E^{2r+1}_{2^{b+1},t+2^{b+1}}$. By assumption, we have
$$
(t+2^{b+1})-2^{b+1}-(2r+1) \leq t+2^{b}+r-2^{b+1}-2r-1 < t-2^{b}-r.
$$
Thus, by induction, $\theta(s,0)x$ is an 
infinite cycle for some $s\gg b+1$.

\bigskip

\noindent ${\bf r=2^{a}-2^{b}:}$ Write $\xi = \sum a_{I}\delta_{I}$ as a 
homogeneous linear combination of admissible operations, with each 
$a_{I}\in \ell$. Then a typical indexing $I$ 
can be written as $I = (i_{1},\ldots,i_{a})$. Admissibility implies 
that $$i_{1}\geq 2i_{2}\geq \ldots 2^{a-1}i_{a}\geq 2^{a}>2^{b}.$$ 

\prop{annprop} now applies to tell us that $\theta(s,b)\xi = 0$ for some 
$s\gg 0$. Thus $\theta(s,0)x=\theta(s,b)\theta(b,0)x$ survives to 
$E^{2^{e}r}$, where $e=s-b$, and $$d_{2^{e}r}\lb\theta(s,b)\theta(b,0) 
x \rb = \theta(s,b) \lb d_{r} \theta(b,0) x \rb = 0 \in E^{2^{e}r}.$$ Thus 
$\theta(s,0)x = \theta(s,b)\theta(b,0)x$ survives to 
$E^{2^{e}r+1}$, so proceed as per the previous case.
\bigskip

Let $y\in \pi_{*}A$ be the element determined by $\lb \theta(s,0) x\rb\in 
E^{\infty}A$. To see that it is non-trivial, suppose that $\lb\theta(s,0) 
x\rb\in 
E^{r}A$ is a boundary. Then $E^{r}(f)(\lb\theta(s,0)x\rb) = 
\lb\theta(s,0) x\rb$ 
is also a boundary in $E^{r}B$. But, since $s\geq b$, this contradicts 
the induction hypothesis.

We conclude, by induction and \eqn{indstep}, that
$$
\Phi_{A}(y) = \Phi_{B}(\alpha_{n-2}^{s-b}(y^{\prime})) = 
\alpha_{n-2}^{s-b}\Phi_{B}(y^{\prime}) = 
\alpha^{s-b}_{n-2}\alpha^{b}(x) = \alpha_{n-2}^{s}(x). 
$$

\hfill $\Box$

\bigskip

\section{Proof of Theorem B}

\begin{theorem}\label{thmb}
	Let $R\to S\to \ell$ be a surjective homomorphisms of Noetherian rings, 
	with $\ell$ a field of characteristic 2, such that $D_{s}(S|R;\ell) = 0$ 
	for $s\gg 0$. Then the following hold:
	\begin{enumerate}
		\item If the divided square $\gamma_{2}$ acts nilpotently on 
		$Q\tor_{*}^{R}(S,\ell)$ it follows that $D_{\geq 2}(S|R;\ell) = 0$.
		
		\item If the Andr\'e operation $\vartheta$ acts nilpotently on 
		$Q\tor^{R}_{*}(S,\ell)$ it follows that $D_{\geq 3}(S|R;\ell)\\ = 0$.
	\end{enumerate}
\end{theorem}

\noindent {\it Proof:}  Using the simplicial model structure for 
simplicial commutative algebras \cite[\S II]{Qui1}, let 
$R\hookrightarrow \Re \stackrel{\sim}{\rightarrow\!\!\!\!\rightarrow} S$ 
be a factorization of $R\to S$ as a cofibration followed by an 
acyclic fibration. Let $A = \Re\otimes_{R}\ell$. Since $\Re$ can be 
chosen to be a degreewise free $R$-algebra, $A$ is a connected 
simplicial supplemented $\ell$-algebra with the properties
\begin{equation}\label{aprop1}
	H^{Q}_{*}(A) \cong D_{*}(S|R;\ell)
\end{equation}
and
\begin{equation}\label{aprop2}
	\pi_{*}A \cong \tor_{*}^{R}(S,\ell).
\end{equation}	
which follows from \cite[(4.7)]{Qui3} and \cite[(4.7)]{Goe}.
The result now follows from Theorem A. \hfill $\Box$

\bigskip

\noindent {\it Proof of Theorem B:}  First, if $D_{\geq 3}(S|R;\ell) = 0$ 
then a Quillen spectral sequence argument shows that 
$\tor_{*}^{R}(S,\ell) \cong \Gamma[D_{*}(S|R;\ell)]$. It follows from 
the Adem relations that $\vartheta$ acts trivially, hence nilpotently, 
on the indecomposables. The converse follows from \thm{thmb}.

Next, $\varphi$ is a complete 
intersection if and only if $D_{2}(S|R;\ell) = 0$, by \cite[VI.25]{And}. Thus  
$\gamma_{2}$ acts trivially, and hence nilpotently, on the 
indecomposables if $\varphi$ is a complete intersection. Conversely, 
if $\gamma_{2}$ acts nilpotently on the indecomposables, then 
$\varphi$ is a complete intersection, by \thm{thmb}. \hfill $\Box$

\hfill 

\section{Proof of Theorem C}

\subsection{A chain condition for the nilpotency of the Andr\'e 
operation}

Let $(R,\fm,\ell) \\ \to S$ be a surjective homomorphism of local rings of 
characteristic 2. Using the simplicial model structure for 
simplicial commutative algebras \cite[\S II]{Qui1}, let 
\begin{equation}\label{cofac}
R\hookrightarrow \Re \stackrel{\sim}{\rightarrow\!\!\!\!\rightarrow} S
\end{equation}
be a factorization of $R\to S$ as a cofibration followed by an 
acyclic fibration. Then the chains $C(\Re)$ provides a free 
resolution of $S$ under $R$. Furthermore, we have
\begin{equation}\label{chainprop}
C(\Re\otimes_{R}\ell) = C(\Re)\otimes_{R}\ell \cong \Frac{C(\Re)}{\fm C(\Re)}
\end{equation}
for which the homology is isomorphic to $\tor^{R}_{*}(S,\ell)$.

We now provide a chain level condition which insures the nilpotency of 
the Andr\'e operation on $\tor$.

\begin{lemma}\label{nilcond}
	If the divided square $\gamma_{2}$ acts nilpotently on $\fm C_{\geq 2}(\Re)$, 
	then the Andr\'e operation $\vartheta$ acts nilpotently on 
	$Q\tor^{R}_{*}(S,\ell)$.
\end{lemma}

\noindent {\it Proof:} Let $x \in C_{\geq 2}(\Re)\otimes_{R}\ell$ be a cycle 
which represents a non-trivial element in homology. Let $y\in C(\Re)$ 
be a pre-image of $x$ under the canonical projection. Since $C(\Re)$ 
is acyclic then $\partial y \in \fm C(\Re)$ by \eqn{chainprop} and is 
non-trivial. 

By \prop{dwylmma} and an induction argument, we have
$$
\partial \alpha_{1}^{s}(y) = \gamma_{2}^{s}(\partial y).
$$
By assumption $\gamma^{s}(\partial y) = 0$ for $s\gg 0$. Thus 
$\alpha_{1}^{s}(y)$ is a boundary for $s\gg 0$ by the acyclicity of 
$C(\Re)$. Hence, by the naturality of $\alpha_{1}$, $\alpha_{1}^{s}(x)$ 
represents the trivial element in the homology of 
$C(\Re\otimes_{R}\ell)$. \hfill $\Box$

\bigskip

\subsection{Some properties of local rings and their homomorphisms}

We now summarize some important properties of certain local rings and 
homomorphisms between them. First, recall that a {\it Cohen-Macaulay ring} is 
a Noetherian ring $B$ such that $\depth B = \dim B$. Cf. \cite{BH}. 

\begin{lemma}\label{locprop}
	\begin{enumerate}
		
		\item  $R$ is a Cohen-Macaulay ring if and only if $R_{\wp}$ and $\widehat 
		R_{\wp}$ are Cohen-Macaulay rings for any $\wp \in \spec R$.
		
		\item If $R$ is a Cohen-Macaulay ring, then any polynomial ring 
		$R[x_{1},\ldots,x_{n}]$ over $R$ is Cohen-Macaulay.
		
		\item If $(R,\fm) \to (S,\fn)$ is a surjective homomorphism of 
		local rings there exists a commutative diagram of local rings
		$$
        \begin{array}{ccc}
         R & \stackrel{\phi}{\rightarrow\!\!\!\!\rightarrow} 
         & (S,\fn) \\[1mm]\downarrow \hspace*{1pt}
         &&
         \hspace*{1pt} \downarrow \\[1mm]
          R^{\prime} & \stackrel{\phi^{\prime}}
         {\rightarrow\!\!\!\!\rightarrow} & (S^{\prime},\fn^{\prime}) \\[1mm]
         \end{array}
         $$
		 with the following properties:
		 \begin{enumerate}
            \item For $s\geq 2$, there exists an injection $D_{s}(S|R;-)\to 
             D_{s}(S^{\prime}|R^{\prime};-)$ which is an isomorphism for $s>2$.
			
			 \item If $S$ is Cohen-Macaulay then $S^{\prime}$ is Artin and 
			 $\depth R^{\prime} = 0$.
			 
			\item $R$ and $S$ are both Cohen-Macaulay if and only if
			$R^{\prime}$ and $S^{\prime}$ are both Artin local rings.
		 \end{enumerate}
	
	 \end{enumerate}
\end{lemma}

\noindent {\it Proof:} For 1., see \cite[2.1.3 and 2.1.8]{BH}. For 2., 
see \cite[2.1.9]{BH}.

For 3., form a surjection $S\to 
(S^{\prime},\fn^{\prime})$ by quotienting out by the ideal generated 
by the maximal regular sequence 
in a minimal generating set for $\fn$. Let $I\subset R$ be the kernel 
of $R\to S\to S^{\prime}$. 
Form a surjection $R\to (R^{\prime},\fm^{\prime})$ by, again, quotienting out 
by an ideal generated by 
the maximal regular sequence in a minimal generating set for $I$. Then 
there is a resulting commuting diagram for which the vertical maps are 
complete intersections (i.e the kernels are generated by regular 
sequences).

Now, applying the Jacobi-Zariski sequence \cite[V.1]{And} to the diagram 
above, we get two long exact sequences
$$
\ldots \to  D_{s}(S|R;-) \to D_{s}(S^{\prime}|R;-) \to 
D_{s}(S^{\prime}|S;-) \to D_{s-1}(S|R;-) \to \ldots
$$
\noindent and
$$
\ldots \to D_{s}(R^{\prime}|R;-) \to D_{s}(S^{\prime}|R;-) \to 
D_{s}(S^{\prime}|R^{\prime};-) \to D_{s-1}(R^{\prime}|R;-) \to \ldots
$$
From \cite[VI.26]{And}, the long exact sequences reduce to two injections
$$
D_{s}(S|R;-) \to D_{s}(S^{\prime}|R;-)
$$
\noindent and
$$
D_{s}(S^{\prime}|R;-) \to D_{s}(S^{\prime}|R^{\prime};-)
$$
for $s\geq 2$. Furthermore, the first map is an isomorphism in the 
same range, as is the second map for $s>2$. Composing the two gives 
the desired map. 

Next, by 1., $S$ is a Cohen-Macaulay ring if 
and only if $S^{\prime}$ is Cohen-Macaulay. Since 
$\depth S^{\prime} = 0$ \cite[1.2]{BH} we have $\dim S^{\prime} = 0$, 
which occurs if and only if $S^{\prime}$ is
Artinian. Cf. \cite[\S 5]{Mat}. Thus $\fn^{\prime}$ is nilpotent as 
an ideal.

Now, if $x \in \fm^{\prime}$ then $\phi^{\prime}(x) \in \fn^{\prime}$ 
satisfies $$\phi^{\prime}(x^{t}) = \phi^{\prime}(x)^{t} = 0 \quad 
\mbox{for} \quad t\gg 0.$$ Thus $x^{t} \in \Ker \phi^{\prime}$ for $t\gg 
0$. From the construction, $\depth (\Ker \phi^{\prime}) = 0$ so we can 
choose $y\neq 0$ in $\Ker \phi^{\prime}$ such that $yx^{t} = 0$ for 
$t\gg 0$. Choose the smallest $t\geq 1$ such that $yx^{t} = 0$. Then 
$u = yx^{t-1} \neq 0$ satisfies $ux = 0$. We conclude that $\depth 
R^{\prime} = 0$. 

Finally, if both $R$ and $S$ are both Cohen-Macaulay then we can 
conclude, from 3. (b), that $S^{\prime}$ is Artin and that $\dim 
R^{\prime} = 0$, hence $R^{\prime}$ is Artin. The converse follows 
from \cite[2.1.3]{BH}.  \hfill $\Box$

\bigskip 

\subsection{Proof of Theorem C}

We are now in a position to prove Theorem C. We first give a result 
which connects the Artinian property on a local ring to the 
nilpotency of the Andr\'e operation on $\tor$.

\begin{theorem}\label{aaconn}
	Let $(R,\fm,\ell) \to S$ be a surjective homomorphism of local rings 
	of characteristic 2. Then $R$ being Artin implies that the Andr\'e 
	operation $\vartheta$ acts nilpotently on $Q\tor^{R}_{*}(S,\ell)$.
\end{theorem}

\noindent {\it Proof:} Factor $R \to S$ as per (\ref{cofac}).  Let $w\in \fm 
C_{\geq 2}(\Re)$. Then 
$$
w = t_{1}x_{1}+\ldots +t_{n}x_{n} \text{\quad with \quad} t_{1},\ldots,t_{n}\in \fm, 
\; x_{1},\ldots,x_{n}\in C_{\geq 2}(\Re).
$$
From the properties of divided squares, we have
$$
\gamma_{2}^{s}(w) = t_{1}^{2^{s}}\gamma_{2}^{s}(x_{1})+\ldots 
t_{n}^{2^{s}}\gamma_{2}^{s}(x_{n}) \; \text{ modulo decomposables}.
$$
Since $R$ is an Artin local ring, then $\fm^{s} = 0$ for $s\gg 0$, by
\cite[2.3]{Mat}, hence
$$
\gamma_{2}^{s}(w) = 0 \; \text{ modulo decomposables,} \quad s\gg 0. 
$$
Since $\gamma_{2}$ kills decomposables in positive degrees, we 
conclude that $\gamma_{2}$ acts nilpotently on $\fm C_{\geq 2}(\Re)$. The 
result now follows from \lma{nilcond}. \hfill $\Box$

\bigskip

\noindent {\it Proof of Theorem C:} By \cite[S.29]{And}, it is enough 
to show that $D_{\geq 3}(S|R;\ell) = 0$, for any residue field $\ell = 
S_{\fn}/\fn S_{\fn}$, $\fn \in \mbox{Spec}\; S$. By the stability of 
Andr\'e-Quillen homology under localization \cite[V.27]{And} 
we may assume that $\varphi: R \to (S,\fn,\ell)$ is a homomorphism 
of local rings. Further, since we are assuming that $\varphi$ is 
essentially of finite type, there is a factorization 
$$
R \stackrel{\tau}{\rarrow} R[X]_{\mathfrak N} = T
\stackrel{\sigma}{\rarrow} S.
$$
\noindent as per (\ref{factor}).
Since $\tau$ is faithfully flat, flat base change \cite[IV.54]{And} 
tells us that
$$
D_{s}(T|R;\ell) \cong D_{s}(\ell[X]|\ell;\ell) = 0 \quad \mbox{for} 
\quad s\geq 1.
$$
Applying the Jacobi-Zariski sequence \cite[V.1]{And}, we conclude that
$$
D_{s}(S|R;\ell) \cong D_{s}(S|T;\ell) \quad \mbox{for} 
\quad s\geq 2.
$$
Note further that $\fd_{T} S = \fd_{R} S$, by a base change spectral sequence argument, 
and if $R$ is Cohen-Macaulay, then $T$ is also Cohen-Macaulay, 
by \lma{locprop}.1 and \ref{locprop}.2. Thus we may assume that $\varphi = \sigma$.

Now suppose $R$ and $S$ are both Cohen-Maclaulay. Then \lma{locprop}.3
allows us to assume that $R$ is an Artin local ring.  Thus 1. follows 
from \thm{aaconn} and Theorem B. 

Finally, if $\fd_{R} S$ is finite, then $Q\tor_{*}^{R}(S,\ell)$ is 
finite and, hence, possesses a nilpotent action of $\gamma_{2}$. 
Thus $\varphi$ is a complete intersection, by Theorem B, giving us 2. 
\hfill $\Box$

\end{document}